\documentclass[a4paper,11pt]{article}

\usepackage[english]{babel}
\usepackage{amssymb,amsmath}

 \usepackage{graphicx}

\def\qed{\hfill $\square$}

\newtheorem{theorem}{Theorem}
\newtheorem{lemma}{Lemma}
\newtheorem{rem}{Remark}

\numberwithin{equation}{section}
\long\def\symbolfootnote[#1]#2{\begingroup\def\thefootnote{\fnsymbol{footnote}}
\footnote[#1]{#2}\endgroup}

\newcommand{\eps}{\varepsilon}

\newcommand{\R}{{\mathbb{R}}}

\begin{document}

\title{
\vspace{0.0in} {\bf\Large  Variables Scaling to Solve a Singular Bifurcation Problem with}\\
 {\bf\Large Applications to Periodically Perturbed Autonomous Systems}\\
 {\bf \small (Dedicated to Prof. R. Johnson on the occasion of his 60th birthday)}}

\author{{\bf\large Mikhail Kamenskii, Oleg Makarenkov, Paolo Nistri}}

\date{}
\maketitle

\symbolfootnote[0]{Mikhail Kamenskii}
\symbolfootnote[0]{Department of Mathematics, Voronezh State University, 394006 Voronezh, Russia}
\symbolfootnote[0]{e-mail: mikhailkamenski@mail.ru\\}

\symbolfootnote[0]{Oleg Makarenkov}
\symbolfootnote[0]{Department of Mathematics, Imperial College London, London, SW7 2AZ,UK}
\symbolfootnote[0]{e-mail: o.makarenkov@imperial.ac.uk\\}

\symbolfootnote[0]{Paolo Nistri}
\symbolfootnote[0]{ Dipartimento di Ingegneria dell'Informazione,
Universit\`a di Siena, 53100 Siena, Italy.}
\symbolfootnote[0]{e-mail: pnistri@dii.unisi.it}

\noindent{\bf Abstract.} {\small By means of a linear scaling of
the variables we convert a singular bifurcation equation in $\R^n$
into an equivalent equation to which the classical implicit
function theorem can be directly applied. This allows to deduce
the existence of a unique branch of solutions as well as a
relevant property of the spectrum of the derivative of the
singular bifurcation equation along the branch. We use these
results to show the existence, uniqueness and the asymptotic
stability of periodic solutions of a $T$-periodically perturbed
autonomous system bifurcating from a $T$-periodic limit cycle of
the autonomous unperturbed system. This problem is classical,
but the novelty of the method proposed is that it allows us to
solve the problem without any reduction of the dimension of the
state space as it is usually done in the literature by means of
the Lyapunov-Schmidt method.}






\section{Introduction}
In Section 2 we consider an abstract bifurcation equation of the
form
\begin{equation}\label{1}
  \Phi(v,\eps):=P(v)+\eps Q(v,\eps)=0
\end{equation}
where $P\in C^2(\R^n, \R^n), Q\in C^1(\R^n\times[0,1], \R^n)$ and,
for $\eps>0$ sufficiently small, we look for the existence of
zeros $v_\eps$ of the map $\Phi$. Here it is assumed the existence
of a $v_0\in \R^n$ such that $P(v_0)=0$ with the matrix $P'(v_0)$
singular. In other words, we deal with an abstract singular
bifurcation problem in $\R^n$ with a small bifurcation parameter
$\eps>0$. Due to the singularity of $P'(v_0)$ it is not possible
to use directly to (\ref{1}) the classical implicit function
theorem to show the existence and uniqueness of a branch
$\{v_\eps\}$, $\eps>0$ small, of solutions of the equation
$\Phi(v,\eps)=0$.

\noindent In this paper, by means of a linear scaling of the
variables $v\in \R^n$ we convert the problem of finding zeros of
(\ref{1}) to the problem of finding zeros of a map $\Psi(w,\eps)$
for which there exists a unique $w_0\in \R^n$ such that
$\Psi(w_0,0)=0$ and  $\Psi'_w(w_0,0)$ is not singular. Therefore,
the new bifurcation equation $\Psi(w,\eps)=0$ can be solved by
means of the classical implicit function theorem to conclude the
existence and uniqueness of a branch of zeros $\{w_\eps\}$, for
$\eps>0$ small. The advantage and the novelty of the approach is
that getting the equation $\Psi(w,\eps)=0$ does not require
solving any implicit equations which is usually done when applying
the Lyapunov-Schmidt reduction approach  (see \cite{ch}, Ch.~2,
\S~4).

 Our bifurcation equation $\Psi(w,\eps)=0$ is, therefore,
formally different from that given by Lyapunov-Schmidt reduction
(see e.g. \cite{l}). That is why we show in Section 3 that
applying our general result to the perturbed autonomous system
\begin{equation}\label{2}
  \dot x=f(x)+\eps g(t,x,\eps).
\end{equation}
where  $f\in C^2(\R^n, \R^n), g\in C^1(\R \times \R^n \times
[0,1], \R^n)$ is $T$-periodic and $\eps>0$ is small, leads to the
same classical Malkin-Loud (or sometimes called Melnikov)
bifurcation function. We end up, therefore, with the statement
that a well known classical result on the existence, uniqueness
and asymptotic stability of a family of $T$-periodic solution of
(\ref{2}) bifurcating from the $T$-periodic limit cycle $x_0$ of
the autonomous system $\dot x=f(x)$ (see Malkin \cite{m}, Loud
\cite{l}, Blekhman \cite{ble}) follows from our bifurcation
theorem, while avoiding the Lyapunov-Schmidt reduction reduces the
analysis significantly.



A first result in this direction has been obtained by the authors
in \cite{enoc08} by means of a version of the implicit function
theorem for directionally continuous functions, see
\cite{bressan}. The idea of using the linear scaling has been,
therefore, already reported at the conference \cite{enoc08}. But
the approach in \cite{enoc08} is based on the employ of
isochronous surfaces of the Poincar\'e map transversally
intersecting the limit cycle $x_0$ that requires a non-trivial
information about smoothness of these surfaces, while the
considerations in this paper rely on very basic facts of analysis
only.


The paper is organized as follows. In Section 2 we first
reduce the abstract singular bifurcation equation (\ref{1}) to an
equivalent non-singular bifurcation equation, then in Theorem 1 we
provide conditions  under which the non-singular problem satisfies
the assumptions of the classical implicit function theorem.
Furthermore, in Theorem 2 we establish a relevant property of the
spectrum of the derivative of the singular bifurcation equation
along the branch which permits to study the asymptotic stability
of the bifurcating zeros. In Section 3, under the standard
assumption that the Malkin's bifurcation function associated to
(\ref{2}) has non-degenerate zeros, the results stated in Section
2 permit to show (Theorem 3) the existence of a parametrized
family of $T$-periodic solutions of (\ref{2}) bifurcating from the
$T$-periodic limit cycle of the unperturbed system as well as
their asymptotic stability.
The main tools to prove Theorem 3 consist in a representation
formula for the Malkin's bifurcation function in terms of the
$T$-periodic perturbation of the autonomous system and of a formula
for its derivative. These formulas are stated in Lemma 2 and Lemma 3
respectively.

\section{Variables scaling to transform a singular bifurcation
problem into a non-singular one}

Consider the function $\Phi: \R^n\times[0,1]\to \R^n$ defined by
\begin{equation}\label{bif}
  \Phi(v,\eps)=P(v)+\eps Q(v,\eps)
\end{equation}
where $P\in C^2(\R^n, \R^n), Q\in C^1(\R^n\times[0,1], \R^n)$ and
$\eps>0$ is a small parameter.

\noindent In this Section, assuming the existence of $v_0\in \R^n$
such that $P(v_0)=0$ with $P'(v_0)$ singular, we provide a method
to show the existence and the uniqueness of the solution $v_\eps$
of the equation
$$
\Phi(v,\eps)=0
$$
for $\eps>0$ sufficiently small, without using the usual
Lyapunov-Schmidt reduction approach. To this aim we assume the
existence of a linear projector $\Pi:\R^n\to\R^n$ such that
$\mbox{Im}\,\Pi \bigoplus \mbox{Ker}\,\Pi=\R^n,$ $\mbox{Im}\,\Pi$
and $\mbox{Ker}\,\Pi$ are invariant subspaces under $P'(v_0)$ and
$\Pi P'(v_0)=\Pi Q(v_0,0)=0$.

\noindent Since $P'(v_0)$ is singular we cannot apply the
classical implicit function theorem, see e.g. \cite{impl-book},
to study the existence of connected components of zeros of $\Phi$
emanating from $(v_0,0)$. Observe that, in general, as it is shown
in \cite{l} and \cite{mn1}, there could exist several branches of
zeros of $\Phi$ emanating from $(v_0,0)$. In this paper we provide
conditions (which are apparently generic when applying the result
to differential equations, see Section 3) under which the branch
is unique. In particular in Section 3, such conditions are
expressed in terms of the Malkin bifurcation function associated
to (\ref{2}), see \cite{m}. More precisely, in Section 3 we have
$v_0=x_0(\theta_0)$, where $x_0$ is a one parameter curve of zeros
of $P$ and $\theta_0$ is a non-degenerate simple zero of the
Malkin bifurcation function. The approach to achieve this result
is commonly based on the classical Lyapunov-Schmidt reduction
method. In the infinite dimensional case, see \cite{henry} and
more recently \cite{banach-kmn}.

\noindent In this paper we propose a different approach based on
an equivalent formulation of the problem. More precisely, by means
of a scaling of the variables, we rewrite the problem of finding
zeros of $\Phi(v,\eps)$, for $\eps>0$ small, as a non-singular
bifurcation problem to which apply the classical implicit function
theorem. Namely, we associate to the map $\Phi$ the following
function
\begin{equation}\label{phi}
  \Psi(w,\eps)=\dfrac{1}{\eps}\left(\Phi(v_0+\eps w,\eps)-\Pi\Phi(v_0+\eps
  w,\eps)+\dfrac{1}{\eps}\Pi \Phi(v_0+\eps w,\eps)\right),
\end{equation}
for any $w\in\mathbb{R}^n$ and any $\eps>0$, and we look for zeros
of $\Psi$ branching from some $(w_0,0).$ Indeed, as it is easy to
see, $(v,\eps)\in\mathbb{R}^n\times[0,1]$ is a zero of $\Phi$ if
and only if $\left(\dfrac{v-v_0}{\eps},\eps\right)$ is a zero of
$\Psi.$

\vskip0.2truecm \noindent In the sequel the vector space of linear
operators $L:\R^n\to\R^n$ will be denoted by
$\mathcal{L}(\mathbb{R}^n)$. Next Lemma provides the main
properties of the function $\Psi$.


\begin{lemma} \label{lem1} Assume that $P\in
C^2(\mathbb{R}^n,\mathbb{R}^n)$ and  $Q\in
C^1(\mathbb{R}^n\times[0,1], \mathbb{R}^n)$. Let
$v_0\in\mathbb{R}^n$ be such that $P(v_0)=0$ and $P'(v_0)$
singular. Let $\Pi:\mathbb{R}^n\to\mathbb{R}^n$ be a linear
projector invariant with respect to $P'(v_0)$ such that $\Pi
P'(v_0)=\Pi Q(v_0,0)=0.$ Define $\Psi(w,0)$ as follows
\begin{equation}\label{psi}
  \Psi(w,0)=\dfrac{1}{2}
  \Pi P''(v_0)ww+\Pi Q'_v(v_0,0)w+\Pi
  Q'_\eps(v_0,0)+(I-\Pi)P'(v_0)w+(I-\Pi)Q(v_0,0)
\end{equation}
with
\begin{equation}\label{psiprime}
\Psi'_w(w,0)=
  \Pi P''(v_0)w+\Pi Q'_v(v_0,0)+(I-\Pi)P'(v_0).
\end{equation}
Then $\Psi\in C^0(\mathbb{R}^n\times\mathbb{R},\mathbb{R}^n)$ and
$\Psi'_w\in
C^0(\mathbb{R}^n\times\mathbb{R},\mathcal{L}(\mathbb{R}^n)).$
\end{lemma}

\noindent{\bf Proof.} From (\ref{phi}) the Taylor expansion with
the rest in the Lagrange's form leads to
\begin{eqnarray*}
   \Pi \Psi(w,\eps)&=&\frac{1}{\eps^2}\Pi \Phi(v_0+\eps
   w,\eps)=\frac{1}{\eps^2}\Pi(P(v_0+\eps w)+\eps Q(v_0+\eps
   w,\eps))=\\
   &=&\frac{1}{\eps^2}\Pi\left(P(v_0)+\eps P'(v_0)w+\frac{1}{2}\eps^2 P''\left(v_0+\widehat{\eps}(w,\eps)w\right)ww\
     +\eps Q(v_0,0)+\right.\\
     & &+\left.\eps^2
     Q'_v\left(v_0+\widetilde{\eps}(w,\eps)w,\widetilde{\eps}(w,\eps)\right)w+\eps^2
     Q'_\eps
     \left(v_0+\widetilde{\eps}(w,\eps)w,\widetilde{\eps}(w,\eps)\right)\right)
\end{eqnarray*}
and
\begin{eqnarray*}
   (I-\Pi)\Psi(w,\eps)&=&\frac{1}{\eps}(I-\Pi)(P(v_0+\eps w)+\eps
   Q(v_0+\eps w,\eps))=\\
   &=&\frac{1}{\eps}(I-\Pi)\left(P(v_0)+\eps P'\left(v_0+\overline{\eps}(w,\eps)w\right)w+\eps Q(v_0+\eps
   w,\eps)\right),
\end{eqnarray*}
where
$\widehat{\eps}(w,\eps),\widetilde{\eps}(w,\eps),\overline{\eps}(w,\eps)\in[0,\eps].$
Using the fact that $P(v_0)=\Pi P'(v_0)=\Pi Q(v_0,0)=0$ we get
\begin{eqnarray*}
   \Psi(w,\eps)&=&\frac{1}{2}\Pi P''\left(v_0+\widehat{\eps}(w,\eps)w\right)ww
     +\Pi
     Q'_v\left(v_0+\widetilde{\eps}(w,\eps)w,\widetilde{\eps}(w,\eps)\right)w+\\
     & & +
     \Pi
     Q'_\eps
     \left(v_0+\widetilde{\eps}(w,\eps)w,\widetilde{\eps}(w,\eps)\right)+
     (I-\Pi)P'\left(v_0+\overline{\eps}(w,\eps)w\right)w+(I-\Pi)Q(v_0+\eps
   w,\eps).
\end{eqnarray*}
From this formula we conclude that $\Psi\in
C^0(\mathbb{R}^n\times\mathbb{R},\mathbb{R}^n).$

\noindent Let us now prove that $\Psi'_w\in
C^0(\mathbb{R}^n\times\mathbb{R},\mathcal{L}(\mathbb{R}^n)).$ The
Taylor expansion applied to $P'(v_0+\eps w)$ permits to write
\begin{eqnarray*}
  \Pi \Psi'_w(w,\eps)&=& \frac{1}{\eps^2} \Pi(\eps P'(v_0+\eps
w) +\eps^2 Q'_v(v_0+\eps
   w,\eps))= \\&=&\frac{1}{\eps^2}\Pi\left(\eps P'(v_0)+\eps^2 P''(v_0+\widetilde{\eps}(w,\eps)w)w+
    \eps^2 Q'_v(v_0+\eps w,\eps)\right),\\
   (I-\Pi)\Psi'_w(w,\eps)&=& \frac{1}{\eps}(I-\Pi)\left(\eps P'(v_0+\eps
    w)+\eps^2 Q'_v(v_0+\eps w,\eps)\right),
\end{eqnarray*}
where $\widetilde{\eps}(w,\eps)\in[0,\eps].$ Taking into account
that $\Pi P'(v_0)=0$ we have
$$
  \Psi'_w(w,\eps)=\Pi P''(v_0+\widetilde{\eps}(w,\eps)w)w+\Pi
  Q'_v(v_0+\eps w,\eps)+(I-\Pi)P'(v_0+\eps
  w)+\eps(I-\Pi)Q'_v(v_0+\eps w,\eps)
$$ and so
$\Psi'_w(w,\eps)\to \Psi'_w(w_0,0)$ as $w\to w_0$ and $\eps\to 0.$
This concludes the proof.
\qed

\begin{rem} An example of linear projector which is invariant with respect to
$P'(v_0)$ is the Riesz projector $\Pi_R:\R^n\to\R^n$ given by
$$
\Pi_R:=\dfrac{1}{2\pi i}\int_{\Gamma} (\lambda I- P'(v_0))^{-1}\,
d\lambda,
$$
where $\Gamma$ is a circumference centered at $0$ and containing in its interior the
only zero eigenvalue of $P'(v_0)$. In fact, by the Riesz
decomposition theorem the subspaces $\mbox{Im}\,\Pi_R$ and
$\mbox{Ker}\,\Pi_R$ are invariant with respect to $P'(v_0),$
$\mbox{Im}\,\Pi_R \bigoplus \mbox{Ker}\,\Pi_R=\R^n$ and $\Pi_R
P'(v_0)=0$.
\end{rem}
We can now prove the following.

\begin{theorem}\label{th1} Assume that $P\in C^2(\mathbb{R}^n,\mathbb{R}^n)$
and $Q\in C^1(\mathbb{R}^n\times[0,1], \mathbb{R}^n)$. Let
$v_0\in\mathbb{R}^n$ be such that $P(v_0)=0$ and $P'(v_0)$ is
singular. Let $\Pi:\mathbb{R}^n\to\mathbb{R}^n$ be a linear
projector (not necessary one-dimensional) invariant with respect
to $P'(v_0)$ with $P'(v_0)$ invertible on $(I-\Pi)\mathbb{R}^n.$
Finally, assume that $\Pi Q(v_0,0)=0,$  $\Pi P''(v_0)\Pi\,r\;
\Pi\, s=0 \quad \mbox{for any} \quad r,s\in\mathbb{R}^n,$ and that
\begin{equation}\label{INV}
-\Pi
P''(v_0)(I-\Pi)\left(P'(v_0)|_{(I-\Pi)\mathbb{R}^n}\right)^{-1}Q(v_0,0)+\Pi
Q'_v(v_0,0)
\end{equation}
is invertible on $\Pi\mathbb{R}^n.$ Then there exists a unique
$w_0\in\mathbb{R}^n$ such that $\Psi(w_0,0)=0$ and
$\Psi'_w(w_0,0)$ is non-singular.
\end{theorem}

\noindent{\bf Proof.} We start by showing the existence of a
$w_0\in\mathbb{R}^n$ such that $\Psi(w_0,0)=0.$ First, observe
that applying $(I-\Pi)$ to (\ref{psi}) we obtain the map $w \to
(I-\Pi) P'(v_0)w + (I-\Pi) Q(v_0,0)$ and the equation
\begin{equation}\label{eqn1}
  (I-\Pi)P'(v_0)w +(I-\Pi)Q(v_0,0)=(I-\Pi)P'(v_0)(I-\Pi)w +(I-\Pi)Q(v_0,0)=0
\end{equation}
is solvable with respect to $(I-\Pi)w$; in fact by our assumptions
$$
  w_1=-\left(\left.P'(v_0)\right|_{(I-\Pi)\mathbb{R}^n}\right)^{-1}Q(v_0,0).
$$
is the solution of (\ref{eqn1}) with $w_1\in (I-\Pi)\R^n$. Now, we
solve the equation
\begin{equation}\label{eqn2}
  \frac{1}{2}\Pi P''(v_0)(\Pi w+w_1)(\Pi
  w+w_1)+\Pi Q'_v(v_0,0)(\Pi w+w_1)+\Pi
  Q'_\eps(v_0,0)=0
\end{equation}
with respect to $\Pi w.$ By assumption $\Pi P''(v_0)\Pi\,r\; \Pi\,
s=0 \;\; \mbox{for any} \;\; r,s\in\mathbb{R}^n,$ moreover
$P''(v_0)ab=P''(v_0)ba,$ hence we can rewrite equation
(\ref{eqn2}) as follows
$$
  \begin{array}{l}
    \Pi P''(v_0)w_1\Pi w+\Pi Q'_v(v_0,0)\Pi w=
   -\dfrac{1}{2}\Pi P''(v_0)w_1\,w_1 -\Pi
    Q'_v(v_0,0))w_1-\Pi Q'_\eps(v_0,0).
  \end{array}
$$
Since by assumption the operator $\Pi P''(v_0)w_1+\Pi Q'(v_0,0)$
is invertible, the last equation has a unique solution $ w_2$ with
$w_2 \in \Pi\,\R^n$. Hence $w_0=w_2 + w_1$ is a zero of
$\Psi(w,0)$.

\noindent From Lemma~\ref{lem1} we have that $\Psi$ is continuous at
$(w_0,0),$ $\Psi'_w$ exists and is continuous at $(w_0,0).$ To
apply the classical implicit function theorem it remains to show
that $\Psi'_w(w_0,0)$ is non-singular. We argue by contradiction
assuming that there exists $h\not= 0$ such that
\begin{equation}\label{ST1}
\Psi'_w(w_0,0)h=\Pi P''(v_0)w_0 h+\Pi
Q'_v(v_0,0)h+(I-\Pi)P'(v_0)h=0.
\end{equation}
Applying $(I-\Pi)$ to (\ref{ST1}) we obtain $(I-\Pi)P'(v_0)h=0$
that is $(I-\Pi)h=0$ and so $h=\Pi h.$ Therefore,
$$
  \begin{array}{l}
     \Pi P''(v_0)w_0h=\Pi P''(v_0)\Pi w_0 \Pi h+\Pi
     P''(v_0)(I-\Pi) w_0 \Pi h=\\
     \qquad\qquad\quad\;\;\;=-\Pi
     P''(v_0)\left(P'(v_0)|_{(I-\Pi)\mathbb{R}^n}\right)^{-1}Q(v_0,0)\Pi
     h
  \end{array}
$$
and applying $\Pi$ to (\ref{ST1}) we obtain
$$
   -\Pi
   P''(v_0)\left(\left.P'(v_0)\right|_{(I-\Pi)\mathbb{R}^n}\right)^{-1}Q(v_0,0)\Pi
   h+\Pi Q'_v(v_0,0)\Pi h=0.
$$
This contradicts our assumption and the proof is completed.
\qed

\begin{rem} The conclusions of Theorem ~\ref{th1} permit to apply
the classical implicit function theorem to obtain the existence of
a $\delta>0$ such that the equation $\Psi(w,\eps)=0$ has, for any
$\eps\in [0, \delta]$,
 a unique solution $w_\eps$  such that
$\|w_0-w_{\eps}\|\le \delta$. Therefore, for $\eps>0$ small, there
exists  a family $\{w_{\eps}\}$ of zeros of the map $\Psi$ such
that $w_{\eps}\to w_0$ as $\eps\to 0$.

\noindent Moreover, under our regularity assumptions $\eps\to
\Phi'_v(v_0+\eps w_{\eps}, \eps)$ is a continuous map; thus, for
any $\eps>0$ sufficiently small, there exists an eigenvalue
$\lambda_{\eps}$ of $\Phi'_v(v_0+\eps w_{\eps}, \eps)$ with the
property that $\lambda_{\eps}\to 0$ as $\eps\to 0$.
\end{rem}

We are now in the position to formulate the following result.

\begin{theorem}\label{th2}
Assume all the conditions of Theorem~\ref{th1} and that zero is a simple
eigenvalue of $P(v_0)$. Let
$v_0=x(\theta_0)$, where $\theta\to x(\theta)$ is a
$C^2$-parametrized curve of zeros of the map $P$. Let
$\{w_{\eps}\}$ and $\{\lambda_{\eps}\}$ as in Remark 2. Let
$\lambda_*\in\mathbb{R}$ be the eigenvalue of the operator
$\left.\Pi P''(v_0)w_0\right|_{\Pi\mathbb{R}^n}+\left.\Pi
Q'_v(v_0,0)\right|_{\Pi\mathbb{R}^n}.$ Then
  $$\lambda_\eps=\eps \lambda_*+o(\eps).$$
\end{theorem}

\noindent{\bf Proof.} Let $l_\eps$ be the unitary eigenvector of
$\Phi'_v(v_0+\eps w_0,\eps)$ associated to the eigenvalue
$\lambda_\eps,$ namely
\begin{equation}\label{TI}
  \Phi'_v(v_0+\eps w_\eps,\eps)l_\eps=\lambda_\eps l_\eps.
\end{equation}
Clearly,
\begin{equation}\label{clear}
l_\eps\to\dfrac{\dot x_0(\theta_0)}{\left\|\dot
x_0(\theta_0)\right\|}\quad{\rm as}\quad \eps\to 0.
\end{equation}
Now we observe that
$$
  \Psi'_w(w,\eps)=\frac{1}{\eps}\left(\eps\Phi'_v(v_0+\eps
  w,\eps)-\eps\Pi\Phi'_v(v_0+\eps
  w,\eps)+\Pi\Phi'_v(v_0+\eps w,\eps)\right)
$$
and using (\ref{TI}) we get
\begin{equation}\label{FI}
  \Pi\Psi'_w(w_\eps,\eps)l_\eps=\frac{1}{\eps}\Pi\Phi'_v(v_0+\eps
  w_\eps,\eps)l_\eps=\frac{1}{\eps}\lambda_\eps\Pi l_\eps
\end{equation}
for any $\eps>0$ sufficiently small. By Lemma~\ref{lem1} as
$\eps\to 0$ we have
  \begin{eqnarray*}
     \Pi\Psi'_w(w_\eps,\eps)l_\eps&\to&\Pi P''(v_0)w_0\dfrac{\dot
     x_0(\theta_0)}{\left\|\dot x_0(\theta_0)\right\|}+\Pi
     Q'_v(v_0,0)\dfrac{\dot
     x_0(\theta_0)}{\left\|\dot x_0(\theta_0)\right\|}\label{FIbis}.
  \end{eqnarray*}
 \noindent From this, by (\ref{FI}) we have that
$\dfrac{\lambda_\eps}{\eps}\to a\in\mathbb{R}$ as $\eps\to 0$ and
$$
\Pi
     P''(v_0)w_0\dot
     x_0(\theta_0)+\Pi
     Q'_v(v_0,0)\dot
     x_0(\theta_0)=a \dot x_0(\theta_0).
$$
Therefore, $a=\lambda_*,$ and the proof is completed.
\qed

\section{An application to periodically perturbed autonomous equations}

In this Section we show that the results of the
previous Section can be straight apply to the problem of bifurcation 
of asymptotically stable $T$-periodic solutions to $T$-periodically perturbed
autonomous systems. Specifically, by
showing that our function (\ref{INV}) is nothing else than the
Malkin's bifurcation function, as far as periodically perturbed
autonomous systems are concerned, we prove the existence of a
unique branch of asymptotically stable periodic solutions
emanating from the family of periodic solutions represented by
limit cycle $x_0$ of the unperturbed system.

\noindent The system under consideration is the following
\begin{equation}\label{ps}
  \dot x=f(x)+\eps g(t,x,\eps).
\end{equation}
where  $f\in C^2(\R^n, \R^n),\, g\in C^1(\R\times \R^n \times
[0,1], \R^n)$ is $T$-periodic and $\eps>0$ is the bifurcation
parameter. We assume that the unique solution of any Cauchy
problem associated to (\ref{ps}) is defined on $[0,T]$.

\noindent We associate to the unperturbed autonomous system
\begin{equation}\label{us}
\dot x=f(x)
\end{equation}
the Malkin's bifurcation function \cite{m}
$$
  M(\theta)=\int_0^T\left<g(t,x_0(t+\theta),0),z_0(t+\theta)\right>dt
$$
where $\left<\cdot, \cdot\right>$ denotes the usual scalar product
in $\R^n$ and $z_0$ is the $T$-periodic function of the adjoint
system
$$
\dot z=-(f'(x_0(t)))^*z
$$
of the linearized system of
$$
 \dot y=f'(x_0(t))y
$$
of autonomous system (\ref{us}). Let $\theta\in[0,T],$  we
define the projector $\Pi:\R^n\to\R^n$ as follows
$$
  \Pi\xi=\dot x_0(\theta)\left<\xi,z_0(\theta)\right>.
$$
 Finally, we convert the problem of finding $T$-periodic
solutions to (\ref{ps}) into the fixed point problem for the
associated Poincar\'e map $\mathcal{P}_{\eps}$ as illustrated in
the following.
We consider the function $x:[0,T]\times\R^n\times[0,1]\to \R^n$
given by
$$
x(t,v,\eps)=x(t)
$$
for all $t\in[0,T],$ where $x(t)$ is the solution of systems
equation (\ref{ps}). The Poincar\'e map for system (\ref{ps}) is
defined by
$$
\mathcal{P}_\eps(v)=x(T,v,\eps).
$$
The functions $P$ and $Q$ of the previous section are
defined as $P(v)=\mathcal{P}_0(v)-v,$
$Q(v,\eps)=\dfrac{\mathcal{P}_\eps(v)-\mathcal{P}_0(v)}{\eps}$
that leads to
$$
  \mathcal{P}_\eps(v)-v=P(v)+\eps Q(v,\eps).
$$
Observe that, since $P(x_0(\theta))=0$ for any $\theta$, we have
that $P'(x_0(\theta))\,\dot x_0(\theta)=0$ and so
$$
(\mathcal{P}_0)'(x_0(\theta))-I=P'(x_0(\theta))
$$
is a singular $n\times n$ matrix for any $\theta\in [0,T]$.

\vskip0.2truecm \noindent With $x_0, z_0, \Pi, P, Q$ as introduced
before we have the following two results. The first one provides a
representation formula for the Malkin's bifurcation function, the
second one a formula for its derivative.

\begin{lemma}\label{lem3}
For any $\theta\in [0,T]$ the limit $Q(v,0):=\lim_{\eps\to
0}Q(v,\eps)$ exists and
$$
M(\theta)=\left<Q(x_0(\theta),0),z_0(\theta)\right>.
$$
Moreover, $Q\in C^1(\mathbb{R}^n\times[0,1],\mathbb{R}^n).$
\end{lemma}

\noindent{\bf Proof.} Differentiating with respect to time one can
see that the function $y(t)=\dfrac{\partial}{\partial \eps}
x(t,x_0(\theta),\eps)$ evaluated at $\eps=0$ solves, for any
$\theta\in [0,T]$, the Cauchy problem
$$
  \dot y=f'(x_0(t+\theta))y+g(t,x_0(t+\theta),0),\quad y(0)=0.
$$
A direct computation shows that
$$
  \frac{d}{dt}\left<y(t),z_0(t+\theta)\right>=\left<g(t,x_0(t+\theta),0),z_0(t+\theta)\right>
$$
and, integrating over the period, yields
$$
  M(\theta)=\left<y(T),z_0(\theta)\right>=\left<Q(x_0(\theta),0),z_0(\theta)\right>.
$$
\qed

\begin{lemma}\label{lem4}
For any $\theta\in [0,T]$ we have
\begin{equation}\label{R}
M'(\theta)=\left<-P''(x_0(\theta))(I-\Pi)\left(\left.P'(x_0(\theta))\right|_{(I-\Pi)\mathbb{R}^n}\right)^{-1}Q(x_0(\theta),0)\dot
   x_0(\theta)+Q'_v(x_0(\theta),0)\dot
   x_0(\theta),z_0(\theta)\right>.
\end{equation}

\end{lemma}

\noindent{\bf Proof.} By Perron's Lemma \cite{perron} we have that
$$
\left<\dot x(\theta), z_0(\theta)\right>=\left<\dot x(0),
z_0(0)\right>
$$
for any $\theta\in [0,T]$. Without loss of generality we may assume
that $\left<\dot x(0), z_0(0)\right>=1.$ As a consequence, by the
definition of the projector $\Pi$, we get
\begin{equation}\label{pi}
\left<\xi, z_0(\theta)\right>=\left<\Pi\,\xi, z_0(\theta)\right>,
\end{equation}
for any $\theta\in [0,T]$. Therefore
$$
\left<P'(x_0(\theta))h,
z_0(\theta)\right>=\left<\Pi\,P'(x_0(\theta))(I-\Pi)h,
z_0(\theta)\right>=0,
$$
for any $\theta\in [0,T]$ and any $h\in \R^n$. Then, by deriving
with respect to $\theta$, we obtain
$$
  \left<P'(x_0(\theta))h,\dot
  z_0(\theta)\right>= \left<-P''(x_0(\theta))\dot
  x_0(\theta)h,z_0(\theta)\right>,
$$
for any $\theta\in [0,T]$ and any $h\in \R^n$. Therefore, we can
rewrite the left hand side of (\ref{R}) with
$(I-\Pi)\left(\left.P'(x_0(\theta))\right|_{(I-\Pi)\mathbb{R}^n}\right)^{-1}Q(x_0(\theta),0)=h$
as follows
$$
  \left<P'(x_0(\theta))(I-\Pi)\left(\left.P'(x_0(\theta))\right|_{(I-\Pi)\mathbb{R}^n}\right)^{-1}Q(x_0(\theta),0),\dot
  z_0(\theta)\right>+\left<Q'_v(x_0(\theta),0)\dot
  x_0(\theta),z_0(\theta)\right>
$$
or equivalently,
$$
  \left<Q(x_0(\theta),0),\dot z_0(\theta)\right>+\left<Q'_v(x_0(\theta),0)\dot
  x_0(\theta),z_0(\theta)\right>,
$$
which is the derivative of $M(\theta)$ at any $\theta\in [0,T]$
according to the formula given by Lemma~\ref{lem3}.
\qed

\vskip0.2truecm Finally, we can prove the following.

\begin{theorem}\label{th3}
 Assume that there exists $\theta_0\in [0,T]$ such that
 $(\mathcal{P}_0)'(x_0(\theta_0))$ has $n-1$ eigenvalues with
 negative real parts, $M(\theta_0)=0$ and
 $M'(\theta_0)<0.$ Then, for $\eps>0$ sufficiently small,
 equation (\ref{ps}) has a unique
 $T$-periodic solution $x_\eps$ such that $x_\eps(t)\to
 x_0(t+\theta_0)$ as $\eps\to 0$ uniformly in $[0,T]$.
 Moreover the solutions $\{x_\eps\}$ are
 asymptotically stable.
\end{theorem}

\noindent{\bf Proof.} Let $v_0=x_0(\theta_0)$, from
Lemma~\ref{lem3} we have
$$
  \Pi Q(x_0(v_0),0)=\dot
  x_0(\theta_0)\left<Q(v_0,0),z_0(\theta_0)\right>=\dot
  x_0(\theta_0)M(\theta_0)=0.
$$
By (\ref{pi}) we obtain
$$
M'(\theta_0)=\left<-\Pi
P''(v_0)(I-\Pi)\left(\left.P'(v_0)\right|_{(I-\Pi)\mathbb{R}^n}\right)^{-1}Q(v_0,0)\dot
x_0(\theta_0) + \Pi Q'_v(v_0,0)\dot x_0(\theta_0),
z_0(\theta_0)\right>\not= 0,
$$
and so (\ref{INV}) is invertible on $\Pi \R^n$. Moreover, from the
fact that $P(x_0(\theta))=0$ for any $\theta\in [0,T]$, we obtain
that
$$
P''(v_0)\dot x_0(\theta_0)\dot x_0(\theta_0) + P'(v_0)
x''_0(\theta_0)=0
$$
Since $\Pi P'(v_0) x''_0(\theta_0)= \Pi P'(v_0) \Pi
x''_0(\theta_0)=0$ we have that  $\Pi P''(v_0)\;\Pi\,r\; \Pi\,s=0$
for any $r,s\in \R^n$. Therefore, all the conditions of
Theorem~\ref{th1} are satisfied and so, compare Remark 2, equation
(\ref{ps}) has a unique $T$-periodic solution $x_\eps$ satisfying
$$
  \left\|w_0-\frac{x_\eps(0)-v_0}{\eps}\right\|\le\delta,
$$
with $\Psi(w_0,0)=0$. Moreover
$$
\Pi P''(v_0) w_0 \dot x_0(\theta_0) + \Pi Q'_v(v_0,0) \dot
 x_0(\theta_0)= \lambda_*\, \dot x_0(\theta_0).
$$
But
$$
\mbox{sign}\,\lambda_*=\mbox{sign}\,\left<\Pi P''(v_0) w_0 \dot
 x_0(\theta_0) + \Pi Q'_v(v_0,0)
\dot
x_0(\theta_0),z_0(\theta_0)\right>=\mbox{sign}\,M'(\theta_0)=-1
$$
Therefore, from Theorem~\ref{th2} there exists $\lambda_\eps=\eps
\lambda_* +o(\eps)$ eigenvalue of
$(\mathcal{P}_\eps)'(x_\eps(0))-I$. This implies that
$$
   {\rm det}\left((\mathcal{P}_\eps)'(x_\eps(0))-I-\lambda_\eps
   I\right)=0.
$$
Hence, $\rho_\eps=1+\lambda_\eps=1+\lambda_*\eps+o(\eps)$ is an
eigenvalue of $(\mathcal{P}_\eps)'(x_\eps(0))$ converging to 1 as
$\eps\to 0.$ Since $\lambda_*<0,$ then $|\rho_\eps|<1$ for
$\eps>0$ sufficiently small. This ends the proof.
\qed

\vskip0.4truecm\noindent {\bf Acknowledgments.} The first author
acknowledge the support by RFBR 09-01-92429 and 06-01-72552. The
second author is supported by the Grant BF6M10 of Russian
Federation Ministry of Education and U.S. CRDF (BRHE), by RFBR
Grant 09-01-00468, by the President of Russian Federation Young
PhD Student grant MK-1620.2008.1 and by Marie Curie grant
PIIF-GA-2008-221331.
The third author is supported by
INdAM-GNAMPA.

\noindent Finally, we would like to acknowledge that the simple
proof of Lemma~2 was suggested by Rafael Ortega during personal 
communications and it is taken from an his own unpublished manuscript.


\begin{thebibliography}{99}



\bibitem{ble}  Blekhman, I.~I. (1971).
Synchronization of dynamical systems, Izdat. Nauka. Moscow

\bibitem{bressan}  Bressan, A. (1988). Directionally continuous selections
and differential inclusions. Funkcial. Ekvac. {\bf 31}, 459-470.

\bibitem{ch} Chow S. N. and Hale J. K. (1982). Methods of bifurcation theory.
Grundlehren der Mathematischen Wissenschaften. {\bf 251},
Springer-Verlag, New York-Berlin.


\bibitem{hale}  Hale, J.K. (1978). Lyapunov-Schmidt method in differential
equations. In Proceedings of the Tenth Brazilian Colloquium. Poços de
Caldas, 1975. Vol. II, pp. 589-603.

\bibitem{henry} Henry, D. (1981). Geometric theory of nonlinear parabolic
equations, Lecture Notes in Mathematics, {\bf 840},
Springer-Verlag, Berlin-New York.



\bibitem{enoc08}
Kamenskii M., Makarenkov O. and Nistri P. (2008). State variables
scaling to solve the Malkin's problem on periodic oscillations in
perturbed autonomous systems, 6th European Nonlinear Dynamics
Conference, ENOC 2008, June 30-July 4, 2008. Saint Petersburg,
Russia. (http://lib.physcon.ru).

\bibitem{banach-kmn}
Kamenskii M., Makarenkov O. and Nistri P. (2008). Periodic bifurcation
from families of periodic solutions for semilinear differential
equations with Lipschitzian perturbation in Banach spaces. Adv.
Nonlinear Stud. {\bf 8}, 271-288.

\bibitem{impl-book}  Krantz S.G. and  Parks H.R. (2003). The Implicit
Function Theorem, History, Theory and Applications, Birkauser
Boston.

\bibitem{l}
Loud W.~S. (1959). Periodic solutions of a perturbed autonomous
system. Ann. Math.  {\bf 70}, 490-529.


\bibitem{mn1} Makarenkov O. and Nistri P. (2008).
``Periodic solutions for planar autonomous systems with nonsmooth
periodic perturbations. J. Math. Anal. Appl. {\bf 338},
1401-1417.

\bibitem{m}
Malkin I. G. (1949). On Poincar\'e's theory of periodic solutions.
Akad. Nauk SSSR. Prikl. Mat. Meh. {\bf 13}, 633-646. (Russian).

\bibitem{perron}
Perron O. (1930). Die Ordnungszahlen der
Differentialgleichungssysteme. Math. Zeitschr. {\bf 31},
748-766.






\end{thebibliography}
\end{document}